# STOCHASTIC PROCESSES WITH SHORT MEMORY

Dmitry Zhabin[1]

*Department of Higher Mathematics and Mathematical Physics,*
*Tomsk Polytechnical University, Tomsk, Lenina avenue 30, 634004 Russia*

The mathematical model of a linear system with the short memory about own stochastic behavior is proposed. It is assumed that the system is under a continual influence of independent stochastic impulses. In a short memory approximation the expression of the stochastic process is found. An application of the model proposed to capital market processes is examined. The approach allows form a stochastic differential for processes concerned. The analog of the Black-Scholes equation for assets dealt on a market with the memory is expressed.

*JEL: C60, D40, G10*

*Keywords: Stochastic processes, autoregressive models*

---

[1] *E-mail address: zhabin@mph.phtd.tpu.edu.ru*

## 1. Introduction

The theory of stochastic processes took a wide application for modeling of the reality (i.g. Tijms 1994). But it should be noted that the most developed part of the stochastic processes theory is the theory of Markov process (Rogers &Williams 2000). Note that the Markov process is a stochastic process such that the future process realization depends solely on current (today's) process value and doesn't refer to any previous realizations.

Actually, the reality is too complicated that there are a number of processes such that the future depends strongly on the history passed. For instance, people that take decisions on capital market are sensitive to the history of securities in some nontrivial way. It is easy to imagine that any real market is a time-depended nonlinear system (see Garbe 1996, Mandelbrot 1963, and Peters 1994) so any capital market process should be described in appropriate way. It follows that a capital market process is an autoregressive one.

All these facts cause a wide application of models with variable volatility like *ARCH* (e.g. Robert Engle, 2001), *GARCH* (see Bollerslev, 1986; Bollerslev et al., 1993), *EWACS* (see Lucas 1990) and so on. As it follows from work of Naukonen 2002, at present time there doesn't exist any model with variable volatility that exceeding others. A chose of the model is caused mostly by specificity of the task solved.

In present paper we consider a linear system that is under continual influence of independent stochastic impulses. We assume that an evolution of the linear system depends on the memory effect of own previous deviations from mean. We found the expression for the stochastic process that generated by such linear system and consider an application of the model proposed to capital market processes.

## 2. Formulation of the model

Let us consider a linear system $x(t)$ at time $t$. We assume that starting from moment $t_0$ the system is under influence of stochastic impulses $y(t)$. Thus we have

$$x(t) = \int_{t_0}^{t} y(s)ds .$$ (1)

Now suppose that impulses have form $\Delta\eta(t)\delta(t - t_k)$, where $\delta(t)$ is Dirac delta function; $\Delta\eta(t)$ are independent stochastic variables with following mean and variance:

$$\mathrm{M}\{\Delta\eta(t)\} = a_k \Delta t_k , \qquad \mathrm{D}\{\Delta\eta(t)\} = b_k^2 \Delta t_k .$$ (2)

Here $a_k$ and $b_k$ are numerical sequences such that all $b_k > 0$. It is clear that linear system (1) could be described in term of stochastic process

$$\xi(t) = \sum_{t_0 \le t_k \le t} \Delta\eta(t_k) .$$ (3)

By assumption that time interval between impulses tents to zero we get the linear system that is under continual influence of independent infinitesimal perturbation. Then we obtain:

$$\xi(t) = \int_{t_0}^{t} d\eta(s) ,$$ (4)

where $d\eta(t)$ are independent stochastic variables with following mean and variance

$$\mathrm{M}\{d\eta(t)\} = a(t)dt , \qquad \mathrm{D}\{d\eta(t)\} = b^2(t)dt .$$ (5)

Hence stochastic process (4) is characterized by mean and variance

$$M\{\xi(t)\} = \int_{t_0}^{t} a(s)ds, \qquad D\{\xi(t)\} = \int_{t_0}^{t} b^2(s)ds.$$ (6)

Without loss of generality it can be assumed that process (4) can be written in form:

$$\xi(t) = \int_{t_0}^{t} a(s)ds + \int_{t_0}^{t} b(s)dw(s),$$ (7)

where the stochastic process $dw(t)$ is defined in the following way

$$dw(t) = \frac{d\eta(t) - M\{d\eta(t)\}}{b(t)}.$$

The integrand of the first item in right side of expression (7) we will name as deterministic item. As above, the integrand of the second item in right side of expression (7) we will name as stochastic item.

Further we assume that the dynamic of linear system (1) is sensitive to deviations of stochastic process (7) from own mean at moments passed. Such effect we can call a "memory effect" (in this case we will tell about process's memory as well). We can assume without loss of generality that such kind of memory depends on value of deviations. Moreover, there exist a specific time $\tau$ that defines a "size" or "depth" of memory to the past. The influence of deviations can be described with help of weight function $f(t-s,\tau)$, where $t$ is the current time, $s$ is some time passed ($s < t$). In the general case weight function is defined from specificity of the tasks solving. But we can write down some general conditions for this function:

$$\begin{cases} f(0,\tau) = 1 \\ f(-\infty,\tau) = 0 \\ \lim_{\tau \to 0+} f(t-s,\tau) = 0 \end{cases} \quad (8)$$

The total contribution of memory effect is the time-averaged value that can be accounted as additive term to stochastic process (7). So we get

$$\xi(t) = \int_{t_0}^{t} a(s)ds + \int_{t_0}^{t} b(s)dw(s) + \frac{1}{t-t_0}\int_{t_0}^{t} f(t-s,\tau)\big(\xi(s) - M\{\xi(s)\}\big)ds . \quad (9)$$

The reader will have no difficulty in showing that the value of memory effect decreases as well as observation period $(t - t_0)$ increases. It follows easily that the mean of process (9) is equal to

$$M\{\xi(t)\} = \int_{t_0}^{t} a(s)ds . \quad (10)$$

Hence for process (9) we obtain

$$\xi(t) = \int_{t_0}^{t} a(s)ds + \int_{t_0}^{t} b(s)dw(s) + \frac{1}{t-t_0}\int_{t_0}^{t} f(t-s,\tau)\left(\int_{t_0}^{s} b(x)dw(x) + \Phi(\xi,s,\tau)\right)ds , \quad (11)$$

where

$$\Phi(\xi,s,\tau) = \frac{1}{t-t_0}\int_{t_0}^{s} f(t-x,\tau)\big(\xi(x) - M\{\xi(x)\}\big)dx . \quad (12)$$

The term $\dfrac{1}{t-t_0}\int_{t_0}^{t} f(t-s,\tau)\Phi(\xi,s,\tau)ds$ describes an influence of memory effects of a second order, i.e., describes the deviations from mean that are generated by previous memory effect. Taking into account that the memory size $\tau$ is assumed to be less then observation time, we can suppose that the contribution of second order memory effect is negligible. It follows that for the linear system with sort memory size we obtain:

$$\xi(t) = \int_{t_0}^{t} a(s)ds + \int_{t_0}^{t} b(s)dw(s) + \frac{1}{t-t_0}\int_{t_0}^{t} ds \int_{t_0}^{s} dw(x)\, b(x) f(t-s,\tau)\,. \qquad (13)$$

By changing the order of integration in (13) we get

$$\xi(t) = \int_{t_0}^{t} a(s)ds + \int_{t_0}^{t} b(s)\left(1 + \frac{1}{t-t_0}\int_{s}^{t} f(t-x,\tau)\,dx\right)dw(s)\,. \qquad (14)$$

### 3. Processes of a capital market

It is often assumed that asset prices must move randomly because of efficient market hypothesis. Let us remind that the most commonly used mathematical representation of asset price is a presentation in term of diffusive processes. Hence, we can write down a stochastic differential for asset price (see Black F. and Scholes 1973, Hull 1996, Wilmott 1993):

$$\frac{\delta S(t)}{S(t)} = a(t)\delta t + b(t)\delta W(t)\,. \qquad (15)$$

Here $a(t)$ and $b(t)$ are piecewise-smooth functions, $\delta W(t)$ is a Wiener process.

It should be noted that the presentation of a market's process as a linear system that is under continual influence of independent stochastic impulses doesn't take into account a possible sensitivities to the history of such processes. But it is easy to imagine that this memory takes place because of people that deal on the market keep in mind previous asset prices. A wide use of the *ARCH* and similar to *ARCH* models also confirm this assumption.

Without loss of generality we can assume that the random walk of the asset price is caused by normal random change of external and unexpected information. This means that asset price can be mathematically presented as

$$S(t) = S(t_0)\exp(h(t)).$$ (16)

Here $h(t)$ – some stochastic process which deterministic and stochastic terms we will note by $A(t)$ and $b(t)$ accordingly. By previous statement, we will assume that the value of each realization of the process $h(t)$ at time $t$ is defined by sum of deterministic, stochastic terms and by the additive term that describe the possible memory about previously passed deviation from the mean. So, we get

$$h(t) = \int_{t_0}^{t} A(s)ds + \int_{t_0}^{t} b(s)\left(1 + \frac{1}{t - t_0}\int_{s}^{t} f(t - x, \tau)dx\right)dW(s).$$ (17)

To analyze a dynamic of such capital market process it is useful to rewrite process (17) in term of stochastic differential $\delta h(t) = h(t + \delta t) - h(t)$. Since (17) is the diffusion process the following stochastic differential take place

$$\delta h(t) = A(t)\delta t + \left(b(t) + \frac{1}{t - t_0}\int_{t_0}^{t} b(s)\left(f(t - s, \tau) - \frac{1}{(t - t_0)^2}\int_{s}^{t} f(t - x, \tau)dx\right)ds\right)\delta W(t).$$ (18)

Let us note that expression (18) is writing down as the limit of equidistance time grid when distance between time points tends to zero.

In spite of the fact that it is impossible to correctly predict future states of the market we can evaluate future value of the mean $A(t)$ and volatility $b(t)$ using bonds prices and implied volatility surface (see Wilmott 1993, Hull 1996).

Let $\widetilde{B}(t)$ be the volatility of process (18)

$$\widetilde{B}(t) = b(t) + \frac{1}{t - t_0} \int\limits_{t_0}^{t} b(s) \left( f(t - s, \tau) - \frac{1}{(t - t_0)^2} \int\limits_{s}^{t} f(t - x, \tau)\, dx \right) ds \,. \qquad (19)$$

It is obvious that as $\tau$ becomes small enough the second item in integrand (19) becomes negligible faster than first item of the integrand. So, an asymptotic behavior of volatility (19) for small $\tau$ is given by

$$\widetilde{B}(t) = b(t) + \frac{1}{t - t_0} \int\limits_{t_0}^{t} b(s) f(t - s, \tau)\, ds \,. \qquad (20)$$

Now if we recall (16), we get

$$\delta S(t) = \left( A(t) + \frac{1}{2} \widetilde{B}^2(t) \right) S(t)\, \delta t + \widetilde{B}(t) S(t)\, \delta W(t) \,. \qquad (23)$$

Expression (23) defines a diffusion process with sort memory about stochastic deviation from mean.

Using (23), it is easy to write down the Black-Scholes equation for asset dealt on the market with sort memory effect. Therefore, we have

$$\frac{\partial V}{\partial t} + \frac{1}{2} \widetilde{B}^2(t) S^2 \frac{\partial^2 V}{\partial S^2} + S \frac{\partial V}{\partial S} - rV = 0 \,. \qquad (24)$$

Here $V = V(t, S)$ is an option price, $r$ is a risk-free interest rate.

Let us consider the case when function функция $f(t - s, \tau)$ is given by following expression

$$f(t, s, \tau) = \exp\left( -\frac{(t - s)^2}{\tau^2} \right) \qquad (25)$$

Then for volatility (19) we get

$$\widetilde{B}(t) = b(t) + \frac{1}{t - t_0} \int\limits_{t_0}^{t} b(s) \left( \exp\left( -\frac{(t - s)^2}{\tau^2} \right) - \frac{\tau \sqrt{\pi}}{2(t - t_0)} \mathrm{Erf}\left( \frac{t - s}{\tau} \right) \right) ds \,, \qquad (26)$$

where Erf($x$) is the well known error function.

It can easily be checked that as $\tau$=0 the second item in the expression (26) is equal to zero and the stochastic differential become describe a dynamic of the stochastic process with volatility $b(t)$. The increasing of $\tau$ causes the increasing of volatility of process $h(t)$. This volatility asymptotically tends to volatility of external impulses $b(t)$ for large observation period ($t - t_0$).

## 4. Some conclusions and remarks

In the present paper we showed that any linear system that is under continual influence of independent external impulses and with short memory effect that is caused by previous deviations from mean at time passed can be described with help of stochastic process with variable volatility. Let us note that such kind models are wide used to investigate an autoregression of stochastic time series. In contrast to *ARCH* and similar to *ARCH* models we didn't make any additional suggestion about behavior of stochastic process.

We showed that the method proposed allows describe processes of capital market and is in accordance with no arbitrage policy. Using the obtained result we wrote down the Black-Scholes equation for a option price. It is easy to see that for vanilla option and for weight function $f(t-s,\tau)$ in form (25) the option price increases together with volatility (26).

John Wiley & Sons Inc., 1993.